\documentclass[11pt,reqno]{amsart}
\usepackage{amssymb}

 \newtheorem{thm}{Theorem}[section]
 \newtheorem{cor}[thm]{Corollary}
 \newtheorem{lem}[thm]{Lemma}

 \theoremstyle{definition}
 
 \theoremstyle{remark}
 
 \numberwithin{equation}{section}

\begin{document}

\title[curvature tensors of Ricci flow]
{On the conditions to control curvature tensors of Ricci flow}

\author{Li Ma and Liang Cheng}


\dedicatory{}
\date{}

\thanks{The research is partially supported by the National Natural Science
Foundation of China 10631020 and SRFDP 20060003002}


\address{Li Ma, Department of Mathematical Sciences, Tsinghua University,
 Peking 100084, P. R. China}

\email{lma@math.tsinghua.edu.cn} \maketitle

\begin{abstract}
An important problem in the study of Ricci flow is that what are
the weakest conditions that provide control of the norm of the
full Riemannian curvature tensor? In this paper, supposing
$(M^n,g(t))$ is a solution to the Ricci flow on a Riemmannian
manifold on time interval $[0,T)$, we show that $L^\frac{n+2}{2}$
norm bound of scalar curvature and Weyl tensor can control the
norm of the full Riemannian curvature tensor if $M$ is closed and
$T<\infty$. Next we prove, without condition $T<\infty$, that
$C^0$ bound of scalar curvature and Weyl tensor can control the
norm of the full Riemannian curvature tensor on complete
manifolds. Finally, we show that to the Ricci flow with bounded
curvature on a complete non-compact Riemannian manifold with the
Ricci curvature tensor uniformly bounded by some constant $C$ on
$M^n\times [0,T)$, the curvature tensor stays uniformly bounded on
$M^n\times [0,T)$. Some other results are also presented.

\end{abstract}

\section{Introduction}
 One of important question in the research of Ricci flow
\begin{equation}\label{Ricciflow}
\frac{\partial g_{ij}}{\partial t}=-2 R_{ij},
\end{equation}
is that what are the weakest conditions to control of the norm of
the Riemannian curvature tensor? There are many results about this
issue. In \cite{SN}, by a blow up argument, N.Sesum showed that
Ricci curvature uniformly bounded on $M \times[0,T)$, where
$T<\infty$, is enough to control the norm of Riemannian curvature
tensor on closed manifolds. Then R.Ye (\cite{RY} and \cite{RY1}) and
B.Wang \cite{BW}, by different arguments, proved that an integral
bound $||Rm||_{\frac{n+2}{2}} = (\int^T_0\int_M |Rm|^{\frac{n+2}{2}}
d\mu dt)^{\frac{ 1}{ \alpha}}$ can control the norm of hole
curvature tensor on closed manifolds, where $Rm$ is the Riemannian
curvature tensor and $T<\infty$. Moreover, B.Wang \cite{W} proved
that if the Ricci tensor is uniformly bounded from below on $[0,T)$
where $T<\infty$, and $||R||_{\alpha} = (\int^T_0\int_M |R|^{\alpha}
d\mu dt)^{\frac{ 1}{ \alpha}}, \alpha\geq \frac{n+2}{2}$, then the
norm of Riemannian curvature tensor is bounded.

Recall that in any dimension $n\geq 3$, the Riemannian curvature
tensor admits an orthogonal decomposition
$$
Rm = -\frac{R}{2(n-1)(n-2)}g\odot g + \frac{1}{n-2}Rc\odot g +W,
$$
where $W$ is the Weyl tensor and $\odot$ denotes the Kulkarni-Nomizu
product. Now we denote
$|Rm|=(g^{i_1j_1}g^{i_2j_2}g^{i_3g_3}g^{i_4j_4}Rm_{i_1i_2i_3i_4}Rm_{j_1j_2j_3j_4})^{\frac{1}{2}}$
and
$|W|=(g^{i_1j_1}g^{i_2j_2}g^{i_3g_3}g^{i_4j_4}W_{i_1i_2i_3i_4}W_{j_1j_2j_3j_4})^{\frac{1}{2}}$
be the norm of the Riemannian curvature tensor and Weyl tensor.

Our first theorem shows that $||R||_{\frac{n+2}{2}}$ and
$||W||_{\frac{n+2}{2}}$ instead of $||Rm||_{\frac{n+2}{2}}$ are
enough to control the norm Riemmannian curvature tensor.

\begin{thm}\label{thm1}
Let $(M^n,g(t))$, $t\in[0,T)$, where $T<\infty$, be a solution to
the Ricci flow (\ref{Ricciflow}) on a closed manifold satisfying
$$||R||_{\frac{n+2}{2}} = (\int^T_0\int_M |R|^{\frac{n+2}{2}} d\mu
dt)^{\frac{2}{n+2}} < \infty$$ and $$||W||_{\frac{n+2}{2}} =
(\int^T_0\int_M |W|^{\frac{n+2}{2}} d\mu dt)^{\frac{2}{n+2}} <
\infty. $$ Then $\sup\limits_{M\times[0,T)}|Rm|<\infty$.
\end{thm}

Note that most of results about this issue should need the
conditions that $M^n$ is closed and $T<\infty$. It is because that
there is no corresponding Perelman's non-local collapsing theorem
for infinite time singularities or non-compact manifolds without
injectivity radius bound conditions. Fortunately, to overcome this,
we can use another simple precompactness theorem (see lemma
\ref{precompactness}) instead of using Perelman's non-local
collapsing theorem. Then we have following theorem without
assumptions $M^n$ is closed and $T<\infty$.

\begin{thm}\label{thm2}
Let $(M^n,g(t))$, $t\in[0,T)$, where $T$ may be infinite, be a
solution to the Ricci flow (\ref{Ricciflow}) on a complete manifold
satisfying $\sup\limits_{M\times[0,T)}|R|<\infty$ and
$\sup\limits_{M\times[0,T)}|W|<\infty$, then
$\sup\limits_{M\times[0,T)}|Rm|<\infty$.
\end{thm}
Next we consider the k\"{a}hler Ricci flow,
\begin{equation}\label{kahler}
\frac{\partial g_{i\bar{j}}}{\partial t}=g_{i\bar{j}}-R_{i\bar{j}},
\end{equation}
on a closed manifold $M$ with $c_1(M)>0$. H.D.Cao proved in
\cite{HC} that (\ref{kahler}) has a solution for all time. One of
the most important questions regarding the K\"{a}hler Ricci flow
(\ref{kahler}) is whether it develops singularities at infinity,
that is whether the curvature of $g(t)$ blows up as $t\to\infty$.
Then Perelman (see \cite{NT}) made a surprising result that the
diameter and scalar curvature of g(t) does not blow up as
$t\to\infty $ for (\ref{kahler}). A interesting application to
theorem \ref{thm2} is following corollary.

\begin{cor}\label{cor1}
Let $(M^n,g(t))$ be a solution to the K\"{a}hler Ricci flow
(\ref{kahler}) on a closed k\"{a}hler manifold with $c_1(X)>0$.
Then $\sup\limits_{M\times[0,+\infty)}|Rm|<\infty$ if and only if
$\sup\limits_{M\times[0,+\infty)}|W|<\infty$.
\end{cor}

Finally, we can extend the result of N.Sesum \cite{SN} about compact
Ricci flow to complete non-compact case using the similar method of
 theorem \ref{thm2}.

\begin{thm}\label{main}
Let $(M^n,g(t))$, $t\in[0,T)$, where $T<\infty$, be a maximal
solution to the Ricci flow (\ref{Ricciflow}) on a complete
manifold $(M^n,g(0))$ satisfying
$\sup\limits_{M\times[0,T)}|Ric|<\infty$ with $(M^n, g(0))$ being
a complete noncompact Riemannian manifold with bounded curvature,
then $\sup\limits_{M\times[0,T)}|Rm|<\infty$.
\end{thm}
We remark that the Ricci flat complete and non-compact manifold is
not in the class in Theorem \ref{main}.

\section{preliminary}\label{preliminary}
 To study the singularities of Ricci flow, R.Hamilton \cite{RH2}
proposed the singularity models which dilate about a singularity
based on the rate of blow up of the curvature. Note that the
curvature bound is immediately satisfied for the blow up about a
singularity, but the injectivity radius bound is not. For finite
time singularities on closed manifolds, Perelman's non-local
collapsing theorem provides the injectivity radius estimate
necessary to obtain a noncollapsing limit. But it is no
corresponding results for infinite time singularities or non-compact
manifolds without injectivity radius bound conditions. In order to
handle this problem, one can use the method, first proposed by
K.Fukaya \cite{F} in metric geometry and D.Glickenstein \cite{DG}
generalized Fukaya's result to Ricci flow, which gives a kind of
precompactness theorem of Ricci flow without injectivity radius
estimates. This method has been used by many people, one can see
\cite{CZ} and \cite{BW} for examples. In fact, we just need a trick
which is used in the proofs in \cite{F} and \cite{DG}. We write them
down below for the convenience of readers and one can see more
general results in \cite{F} and \cite{DG}.

First we need a lemma which is elementary in Riemmanian geometry.

\begin{lem}\label{injestimate}
Let $(M^n, g)$ be a complete Riemannian manifold with $|sec| \leq
1$, $p \in M^n$ and the exponential map $exp_p:T_p M^n \to M^n$ such
that $B(o,\pi) \subset T_p M^n$ equipped with metric $exp_p^*g$.
Then the injectivity radius at o in $B(o,\pi)$ has
$\text{inj}(o)>\frac{\pi}{2}$.
\end{lem}
\begin{proof}
By Rauch comparison theorem, the conjugate radius at $o$ is not
less than $\pi$. Then we can pull back the metric of $M^n$ to
$B(o,\pi)$ with exponential map $exp_p$. By the definition of
exponential map, me know that $exp_p$ map any straight line at $o$
to the geodesic emanating from $p$. Since $B(o,\pi)$ equipped with
metric $exp_p^*g$, the straight line at $o$ must be the geodesic
with respect to the metric $exp_p^*g$ in $B(o,\pi)$. Hence we have
any geodesic at $0$ is straight line and never goes back in
$B(o,\pi)$. Then the lemma follows from Klingerberg's injectivity
estimate immediately.
\end{proof}

A simple application of lemma \ref{injestimate} is following
precompactness theorem for Ricci flow.
\begin{lem}\label{precompactness}
 Let $\{(M^n_i,g_i(t),x_i)\}_{i=1}^{\infty}$, $t\in [0,T]$, be a sequence of
  Ricci flow on complete manifolds such that
$\sup\limits_{M_i\times [0,T]}|Rm(g_i(t))|g_i(t)\leq 1$. Let
$\phi_i=exp_{x_i,g_i(0)}$ be the exponential map with respect to
metric $g_i(0)$ and $B(o_i,\frac{\pi}{2})\subset T_{x_i}M_i$
equipped with metric $\widetilde{g}_i(t)\triangleq \phi_i^*g(t)$.
Then $(B(o_i,\frac{\pi}{2}),\widetilde{g}_i(t),o_i)$ subconverges to
a Ricci flow $(B(o,\frac{\pi}{2}),\widetilde{g}(t),o)$ in
$C^{\infty}$ sense, where $B(o,\frac{\pi}{2})\subset \mathbb{R}^n$
equipped with metric $\widetilde{g}(t)$.
\end{lem}
\begin{proof}
Since $\widetilde{g}_i(t)\triangleq \phi_i^*g(t)$, we have
$\sup\limits_{B(o_i,\frac{\pi}{2})\times
[0,T]}|Rm(\widetilde{g}_i(t))|\widetilde{g}_i(t)\leq 1$. By lemma
\ref{injestimate}, we have
$inj(o_i,\widetilde{g}_i(0))>\frac{\pi}{2}$. Then the lemma follows
from the local version of Hamilton's precompactness theorem for the
Ricci flow (see \cite{BC}, p236).
\end{proof}

Finally, we need a lemma to prove theorem \ref{main}.

\begin{lem}\label{estimate2}
Suppose $(M,g(t))$ is the solution to the Ricci flow for $t\in
[0,T)$, where $T<\infty$, such that the Ricci curvature is uniformly
bounded by some constant C. We have for all $\epsilon > 0$ there
exists an $ \delta(\epsilon,C)> 0$ such that if $ t\in
[t_0,t_0+\delta]$ then
\begin{eqnarray}
\text{Vol}_{g(t)}(B_{g(t)}(x,r))
\geq(1-\epsilon)^{\frac{n}{2}}\text{Vol}_{g(t_0)}(B_{g(t_0)}(x,\frac{r}{1+\epsilon})).
\end{eqnarray}
for all $x\in M$, where $\text{Vol}_{g(t)}(B_{g(t)}(x,r))$ is the
volume of ball centered at point $x$ with respect to metric $g(t)$.
\end{lem}
\begin{proof}
Since the Ricci curvature is uniformly bounded, it is not hard to
show (see lemma 17 in \cite{DG}) that all $\epsilon > 0$ there
exists an $\delta(\epsilon,C)
> 0$ such that if $|t - t_0| < \delta$ then
$$
-\epsilon g(t_0) \leq g(t) - g(t_0) \leq \epsilon g(t_0),
$$
and
$$
|d_{g(t)}(x, y) - d_{g(t_0)}(x, y)| \leq \epsilon^{ \frac{1}{ 2}}
d_{g(t_0)}(x, y),
$$
for all $x, y \in M$. For any $x\in M^n$, $r>0$, and $t\in
[t_0,t_0+\delta]$, we have
$$
B_{g(t_0)}(x,\frac{r}{1+\epsilon}) \subset B_{g(t)}(x,r),
$$
and
$$
\text{Vol}_{g(t)}(B_{g(t_0)}(x,\frac{r}{1+\epsilon}))\geq
(1-\epsilon)^{\frac{n}{2}}\text{Vol}_{g(t_0)}(B_{g(t_0)}(x,\frac{r}{1+\epsilon})).
$$
Hence we have
\begin{eqnarray*}
\text{Vol}_{g(t)}(B_{g(t)}(x,r))
\geq(1-\epsilon)^{\frac{n}{2}}\text{Vol}_{g(t_0)}(B_{g(t_0)}(x,\frac{r}{1+\epsilon})).
\end{eqnarray*}
\end{proof}

\section{proof of theorem \ref{thm1}-\ref{main}}

In this section, we give the proofs to  theorem
\ref{thm1}-\ref{main} and corollary \ref{cor1}. First, we use a
simple blow up argument to prove theorem \ref{thm1}.

\textbf{Proof of Theorem \ref{thm1}.} We argue by contradiction.
Assume that the curvature blows up at a finite time $T$ under the
assumptions of the theorem. Then there is a sequence $(x_i, t_i)$
with $\lim\limits_{i\to\infty} t_i = T$ and $Q_i\to\infty$, where
$Q_i\triangleq |Rm| (x_i, t_i) = \max\limits_{ (x,t)\in
M\times[0,t_i]} |Rm|(x, t)$. We define the pointed rescaled
solutions $(M^n,g_i(t),x_i)$, $t\in [-t_i Q_i,0]$, by $g_i(t)
=Q_ig(t_i+Q_i^{-1}t)$. Then we have a sequence of Ricci flow
$\{(M^n, g_i(t),x_i), -1 \leq t\leq 0\}$ satisfying
$$
\sup\limits_{M^n\times(\infty,0]}|Rm_i|_{g_i}\leq 1
$$
and
$$
|Rm_{i}(x_i,0)|_{g_i}=1,
$$
where $Rm_i$ is the curvature tensor with respect to $g_i$. By
Perelman's non-local collapsing theorem, we have uniform lower bound
of injectivity radius at points $(x_i, t_i)$ for the sequence
$\{(M^n, g_i(t),x_i), -1 \leq t\leq 0\}$. So this sequence
subconverges to an ancient Ricci flow solution $( \overline{M} ,
\overline{x}, \overline{g}(t))$. We calculate
\begin{eqnarray*}
\int^0_{-1}\int_{B_{\overline{g}(0)}(\overline{x},1)}
|\overline{R}|^{\frac{n+2}{2}} d\mu dt &\leq&
\lim\limits_{i\to\infty} \int^0_{-1}\int_{B_{g_i(0)}(x_i,1)} |R_i|^{\frac{n+2}{2}} d\mu dt\\
& =& \lim\limits_{i\to\infty} \int^{t_i}_{t_i-(Q_i)^{-1}}
\int_{B_{g(t_i)}(x_i,Q_i^{-\frac{1}{2}})} |R|^{\frac{n+2}{2}} d\mu dt\\
&\leq& \lim\limits_{i\to\infty} \int^{t_i}_{t_i-(Q_i)^{-1}}\int_{M} |R|^{\frac{n+2}{2}} d\mu dt\\
& =& 0.
\end{eqnarray*}
The last equality holds, since $\int^T_0\int_{M} |R|^{\frac{n+2}{2}}
d\mu dt<\infty$. Likewise,
\begin{eqnarray*}
\int^0_{-1}\int_{B_{\overline{g}(0)}(\overline{x},1)}
|\overline{W}|^{\frac{n+2}{2}} d\mu dt &\leq&
\lim\limits_{i\to\infty} \int^0_{-1}\int_{B_{g_i(0)}(x_i,1)} |W_i|^{\frac{n+2}{2}} d\mu dt\\
& =& \lim\limits_{i\to\infty} \int^{t_i}_{t_i-(Q_i)^{-1}}\int_{B_{g(t_i)}(x_i,Q_i^{-\frac{1}{2}})} |W|^{\frac{n+2}{2}} d\mu dt\\
&\leq& \lim\limits_{i\to\infty} \int^{t_i}_{t_i-(Q_i)^{-1}}\int_{M} |W|^{\frac{n+2}{2}} d\mu dt\\
& =& 0.
\end{eqnarray*}
Since $(\overline{M}, \overline{g}(t))$ is a smooth Riemannian
manifold, we have $\overline{R}(x,t)\equiv 0,
\overline{W}(x,t)\equiv 0$ for $t\in [-1,0]$ in
$B_{\overline{g}(0)}(\overline{x},1)$. By the scalar curvature
evolution equation $\frac{\partial \overline{R}}{\partial t}=\Delta
\overline{R}+2|\overline{Ric}|^2$, we have
$\overline{Ric}(x,t)\equiv 0$ for $t\in [-1,0]$ in
$B_{\overline{g}(0)}(\overline{x},1)$. Then $\overline{Rm}\equiv 0$
in $B_{\overline{g}(0)}(\overline{x},1)$. So we get a contradiction
to $|\overline{Rm}|(\overline{x},0)=1$. $\Box$

Next we prove theorem \ref{thm2} by using blow up arguments and
lemma \ref{precompactness}.

\textbf{Proof of Theorem \ref{thm2}.} The proof is by contradiction.
Then there exist sequences $t_i\to T$ and $x_i\in M$, such that
$$Q_i \doteq |Rm|(x_i, t_i) \geq C^{-1} \sup\limits_{ M\times [0,t_i]} |Rm|(x, t)\to\infty,$$
 as $i\to\infty$, where $C$ is some constant bigger than 1.
 We define the pointed rescaled solutions $(M^n,g_i(t),x_i)$,
 $t\in [-t_i Q_i,0]$, by $g_i(t) =Q_ig(t_i+Q_i^{-1}t)$. Then we have
$$\sup\limits_{M\times (-t_i Q_i,0]} |Rm_i|_{g_i}(x, t)\leq C$$
 and
$$|Rm_{i}(x_i,0)|_{g_i}=1,$$
where $Rm_i$ is the curvature tensor with respect to $g_i$. Let
$\phi_i=exp_{x_i,g_i(0)}$ be the exponential map with respect to
metric $g_i(0)$ and $B(o_i,\frac{\pi}{2})\subset T_{x_i}M_i$
equipped with metric $\widetilde{g}_i(t)\triangleq \phi_i^*g_i(t)$.
So we also have
$$\sup\limits_{B(o_i,\frac{\pi}{2})\times (-t_i Q_i,0]} |\widetilde{Rm}_i|_{\widetilde{g}_i}(x, t)\leq C$$
 and
$$|\widetilde{Rm}_{i}(o_i,0)|_{\widetilde{g}_i}=1,$$
where $\widetilde{Rm}_i$ is the curvature tensor with respect to
$\widetilde{g}_i$. Then
$(B(o_i,\frac{\pi}{2}),\widetilde{g}_i(t),o_i)$ subconverges to a
Ricci flow $(B(o,\frac{\pi}{2}),\widetilde{g}(t),o)$ in $C^{\infty}$
sense by lemma \ref{precompactness}. Since we have that $|R(t)|\leq
C$, where $R(t)$ is a scalar curvature with respect to $g(t)$. Then
$\widetilde{R_i} = \frac{R}{Q_i}$ for all $t \in (-t_iQ_i,0]$, we
get that $\widetilde{R_i}(t)\to 0$ as $i\to \infty$, i.e.
$\widetilde{R}(t)=0$ for $t \in(-\infty, 0]$, where
$\widetilde{R}(t)$ is the scalar curvature with respect to metric
$\widetilde{g}(t)$. By the evolution equation for the scalar
curvature $\frac{d}{ dt}\widetilde{R} = \Delta \widetilde{R} +
2|\widetilde{Ric}|^2$, we have $\widetilde{Ric}(t) = 0$ for all
$t\in(-\infty, 0]$, i.e. our solution $\widetilde{g}$ is stationary.
Likewise, $|\widetilde{W_i}| = \frac{|W|}{Q_i}$ for all $t \in
(-t_iQ_i,0]$, we get that $|\widetilde{W}_i(t)|\to 0$ as $i\to
\infty$, i.e. $\widetilde{W}(t)=0$ for $t \in(-\infty, 0]$. Then we
have $\widetilde{Rm}\equiv0$, which contradict to
$|\widetilde{Rm}(o,0)|=1$. $\Box$

To prove corollary \ref{cor1}, we need following theorem due to
G.Perelman.

\begin{thm}\label{perelman}\cite{NT}
Let $(M^n,g(t))$ be a solution to the k\"{a}hler Ricci flow
$\frac{\partial g_{i\bar{j}}}{\partial t}=g_{i\bar{j}}-R_{i\bar{j}}$
on a compact k\"{a}hler manifold with $c_1(X)>0$. There exists a
uniformly constant such that $|R(g(t))|\leq C$,
$\text{diam}(M,g(t))\leq C$ and $|u|_{C^1}\leq C$, where $u$ is the
Ricci potential satisfying
$g_{i\bar{j}}-R_{i\bar{j}}=\partial_i\partial_{\bar{j}}u$.
\end{thm}

Now we can prove corollary \ref{cor1}.

 \textbf{Proof of Corollary\ref{cor1}.} Corollary \ref{cor1} follows
 from theorem \ref{perelman} and the same arguments of the proofs in
 theorem \ref{thm2}. Note that the evolution equation for the scalar
 curvature $\frac{d}{
dt}R = \Delta R + |Ric|^2-R$ under the k\"{a}hler Ricci flow
(\ref{kahler}). So $R \equiv 0$ still implies $Ric \equiv 0$. $\Box$

Finally, we use the similar methods in N.Sesum \cite{SN} and lemma
\ref{precompactness} to prove theorem \ref{main}.

\textbf{Proof of Theorem \ref{main}.} We argue by contradiction.
Assume that the curvature blows up at a finite time $T$ under the
assumptions of the theorem. Then there exist sequences $t_i\to T$
and $x_i\in M$, such that
$$
Q_i \triangleq |Rm|(x_i, t_i) \geq C^{-1} \sup\limits_{ M\times
[0,t_i]} |Rm|(x, t)\to\infty,
$$
 as $i\to\infty$, where $C$ is some constant bigger than 1.
 We define the pointed rescaled solutions $(M^n,g_i(t),x_i)$,
 $t\in [-t_i Q_i,0]$, by $g_i(t) =Q_ig(t_i+Q_i^{-1}t)$. Then we have
 $\sup\limits_{M\times (-t_i Q_i,0]} |Rm|_{g_i}(x, t)\leq C$
  and $|Rm_{g_i}(x_i,0)|=1$. Let $\phi_i=exp_{x_i,g_i(0)}$ be the exponential map with respect to
metric $g_i(0)$ and $B(o_i,\frac{\pi}{2})\subset T_{x_i}M_i$
equipped with metric $\widetilde{g}_i(t)\triangleq \phi_i^*g_i(t)$.
Then $(B(o_i,\frac{\pi}{2}),\widetilde{g}_i(t),o_i)$ subconverges to
a Ricci flow $(B(o,\frac{\pi}{2}),\widetilde{g}(t),o)$ in
$C^{\infty}$ sense by lemma \ref{precompactness}. Since the Ricci
curvatures of our original flow are uniformly bounded, we have that
$|R(t)|\leq C$, where $R(t)$ is a scalar curvature with respect to
$g(t)$. Since $R_i= \frac{R}{Q_i}$ for all $t \in (-t_iQ_i,0]$, we
get that $R_i(t)\to 0$ as $i\to \infty$, i.e. $\widetilde{R}(t)=0$
for $t \in(-\infty, 0]$. By the evolution equation for a scalar
curvature $\frac{d}{ dt}\widetilde{R} = \Delta \widetilde{R} +
2|\widetilde{Ric}|^2$, we have $\widetilde{Ric}(t) = 0$ for all
$t\in(-\infty, 0]$, i.e. our solution $\bar{g}$ is stationary, where
$\widetilde{R}(t)$ is the scalar curvature with respect to metric
$\widetilde{g}(t)$. Therefore it can be extended for all times $t
\in(-\infty,+\infty)$ to an eternal solution. Take sufficient small
$r > 0$ such that $B(o,r)\subset B(o,\frac{\pi}{2})$. Since
$\widetilde{g}$ is the $C^{\infty}$ limit of $\phi^{*}_ig_i$, we
have
$$
\frac{\text{Vol}B(o,r)}{ r^n} = \lim\limits_{i\to\infty}
\frac{\text{Vol}_{\phi^{*}_ig_i(0)}B_{\phi^{*}_ig_i(0)}(o_i,r)}{r^n}
= \lim\limits_{i\to\infty}
\frac{\text{Vol}_{\phi^{*}_ig(t_i)}B_{\phi^{*}_ig(t_i)}(o_i,
rQ_i^{-\frac{1}{2}})}{(rQ_i^{-\frac{1}{2}})^n},
$$
where $\text{Vol}B(o,r)$ is the volume of ball $B(o,r)$ on
$B(o,\frac{\pi}{2})$ with respect to metric $\widetilde{g}$. For any
$\epsilon>0$, let $\delta(\epsilon,C)$ be the constant defined in
lemma \ref{estimate2}. Since $t_i\to T$ and T is finite, there
exists sufficient large $i_0$ such that
$|t_i-t_{i_0}|<\delta(\epsilon,C)$ and
$$
 \frac{\text{Vol}B(o, r)}{ r^n} \geq \frac{\text{Vol}_{\phi^{*}_ig(t_i)}B_{\phi^{*}_ig(t_i)}(o_i,
rQ_i^{-\frac{1}{2}})}{(rQ_i^{-\frac{1}{2}})^n}-\epsilon,
$$
when $i>i_0$. Hence by lemma \ref{estimate2},
we have
$$
 \frac{VolB(o, r)}{ r^n} \geq
 (1-\epsilon)^{\frac{n}{2}}\frac{\text{Vol}_{\phi^{*}_{i_0}g(t_{i_0})}B_{\phi^{*}_{i_0}g(t_{i_0})}
 (o_{i_0},\frac{1}{1+\epsilon} rQ_i^{-\frac{1}{2}})}{
(rQ_i^{-\frac{1}{2}})^n}-\epsilon.
$$
Note that $\text{Vol}B(p,r)$ has the expansion within the
injectivity radius of $p$ (see \cite{GD} theorem 3.98),
$$
\text{Vol}B(p,r)=\omega_n r^n(1-\frac{R(p)}{6(n+2)}r^2+o(r^2)),
$$
for $B(p,r)$ contained in any complete manifold $M$. Since curvature
is bounded at time $t_{i_0}$ and $\text{inj}(o_i)>\frac{\pi}{2}$, we
have $
\frac{\text{Vol}_{\phi^{*}_{i_0}g(t_{i_0})}B_{\phi^{*}_{i_0}g(t_{i_0})}
 (o_{i_0},\frac{1}{1+\epsilon} rQ_i^{-\frac{1}{2}})}{
(\frac{1}{1+\epsilon}rQ_i^{-\frac{1}{2}})^n} $ uniformly converges
to $\omega_n$ as $i\to\infty$. Hence we have
$$
 \frac{\text{Vol}B(o, r)}{ r^n} \geq (1-\epsilon)^{\frac{n}{2}}
 (\frac{1}{1+\epsilon})^n\omega_n-\epsilon.
$$
Let $\epsilon\to 0$, we conclude that
$$
 \frac{\text{Vol}B(o, r)}{ r^n} \geq \omega_n.
$$
Finally by $\widetilde{Ric}(t) = 0$ and Bishop-Gromov volume
comparison, we have metric $\widetilde{g}$ is flat, which
contradicts to the fact $|\widetilde{Rm}|(o,0)=1$. $\Box$

\end{document}